\documentclass[leqno]{birkjour}
\usepackage[utf8]{inputenc}
\usepackage[USenglish]{babel}

\usepackage{booktabs}

\usepackage{bbm}

\numberwithin{equation}{section}

  \newcommand{\Z}{\mathbb{Z}}
  
  \newcommand{\C}{\mathbbm{C}}
  \newcommand{\RP}{\mathbbm{R}\mathrm{P}}
  \newcommand{\R}{\mathbbm{R}}
  \newcommand{\Rschlange}{\underset{\widetilde{ }}{\mathbbm{R}}}
  \newcommand{\KP}{\mathbbm{K}\mathrm{P}}
  \newcommand{\K}{\mathbbm{K}}
  \renewcommand{\H}{\mathbbm{H}}

  \newcommand{\lt}{\left(}
  \newcommand{\rt}{\right)}
  \newcommand{\lb}{\left\lbrace}
  \newcommand{\rb}{\right\rbrace}

  \newcommand{\id}{\operatorname{id}}
  \newcommand{\pr}{\operatorname{pr}}

  \newcommand{\MC}{\operatorname{MC}}
  \newcommand{\MCC}{\operatorname{MCC}}
  \newcommand{\MF}{\operatorname{MF}}
  \newcommand{\Ker}{\operatorname{Ker}}

  \newcommand{\ST}{\ensuremath{\mathit{ST}}}
  \newcommand{\comp}{\, \scriptstyle \circ \displaystyle \,}
  
    \newtheorem{thm}{Theorem}
    
    \newtheorem*{WeckenTheorem}{Wecken theorem in coincidence theory}

  \theoremstyle{definition}

  \theoremstyle{remark}
    
    \newtheorem{exa}[equation]{Example}
    
    \newtheorem*{remNo}{Remark}
  \newtheoremstyle{example}{3pt}{3pt}{}{}{\bfseries}{:}{.5em}{}
  \theoremstyle{example}
    \newtheorem{MyExa}[equation]{Example}

  \newtheoremstyle{importantquestion}
    {3pt}
    {3pt}
    {\itshape}
    {}
    {\itshape}
    {.}
    {.5em}
    {}
  \theoremstyle{importantquestion}
    \newtheorem*{CQ}{CENTRAL QUESTION}

  \newtheoremstyle{unnumberedtheorem}{3pt}{3pt}{\itshape}{}{\itshape}{}{.5em}{}
  \theoremstyle{unnumberedtheorem}
    \newtheorem*{MythmNo}{Theorem}

\newcommand{\MyVert}{\ensuremath{\; \vert \;}}

\begin{document}

\title{Coincidences and secondary Nielsen numbers}
\author{\textbf{Ulrich Koschorke}}	
\address{%
Department Mathematik\\
University of Siegen\\
57068 Siegen, Germany
}
\email{%
koschorke@mathematik.uni-siegen.de
}

\begin{abstract}
Let
$ f_1, f_2 \colon X^m \longrightarrow Y^n $
be maps between smooth connected manifolds of the indicated dimensions
$ \!m\! $
and
$ \!n \!\!\!$.
Can
$ f_1, f_2 $
be deformed by homotopies until they are coincidence free (i.e.
$ f_1(x) \neq f_2(x) $
for all
$ x \in X $)?
The main tool for addressing such a problem is tradionally the (primary) Nielsen number
$ N(f_1, f_2) $.
E.g. when
$ m < 2n - 2 $
the question above has a positive answer precisely if
$ N(f_1, f_2) = 0 $.
However, when
$ m = 2n - 2 $
this can be dramatically wrong, e.g. in the fixed point case when
$ m = n = 2 $.
Also, in a very specific setting the Kervaire invariant appears as a (full) additional obstruction.

In this paper we start exploring a fairly general new approach.
This leads to secondary Nielsen numbers
$ SecN(f_1, f_2) $
which allow us to answer our question e.g. when
$ m = 2n - 2,\ \; n \neq 2 $
is even and
$ Y $
is simply connected.
\end{abstract}

\thanks{Supported in part by DFG (Deutsche Forschungsgemeinschaft)}



\maketitle


\section{Introduction}

Throughout this paper let
\begin{equation*}
  f_1, f_2 \colon X^m \longrightarrow Y^n
\end{equation*}
be (continuous) maps between smooth connected manifolds without boundary, of the indicated dimensions
$m, n \geq 1, \; X $
being compact.

We are interested in those aspects of the coincidence subspace
\begin{equation*}
  C(f_1, f_2) := \left\lbrace x \in X \MyVert f_1(x) = f_2(x) \right\rbrace
\end{equation*}
in
$ X $
which remain unchanged by homotopies of
$ f_1 $
and
$ f_2 $.
These aspects are reflected to a large extend by the \textit{\textbf{m}inimum numbers}
\begin{align*}
  \MC(f_1, f_2) &:= min \left\lbrace \# C(f_1', f_2') \MyVert f_1' \sim f_1, f_2' \sim f_2 \right\rbrace \\
  \intertext{and, better yet,}
  \MCC(f_1, f_2) &:= min \left\lbrace \# \pi_0(C(f_1', f_2')) \MyVert f_1' \sim f_1, f_2' \sim f_2 \right\rbrace
\end{align*}
of \textit{\textbf{c}oincidence points} and of \textit{\textbf{c}oincidence \textbf{c}omponents}, resp., as the maps vary within the given homotopy classes
$ [ f_1 ], [ f_2 ] $.

The principal problem in topological coincidence theory is to determine these minimum numbers, and --in particular-- to decide when they vanish, i.e. when
$ ( f_1, f_2 ) $
is homotopic to a coincidence free pair.
In this case we say that the pair
$ ( f_1, f_2 ) $
is \textbf{loose}.

Generalizing a well known notion from fixed point theory, we introduced (in \cite{ko2}) a Nielsen number
$ N(f_1, f_2) $
which depends only on the homotopy classes of
$ f_1 $
and
$ f_2 $
and satisfies
\begin{equation*}
  0 \leq N(f_1, f_2) \leq \MCC(f_1, f_2) \leq \MC(f_1, f_2) \leq \infty \;.
\end{equation*}
Furthermore, we proved the following
\begin{WeckenTheorem}\label{MyWeckenTheorem}
  Assume
  $ m < 2n -2 $.\\
  Then for all maps
  $ f_1, f_2 \colon X^m \longrightarrow Y^n $
  we have
  $ \MCC(f_1, f_2) = N(f_1, f_2) \!$.
  In particular,
  $ (f_1, f_2) $
  is loose if and only if
  $ N(f_1, f_2) = 0 $.
\end{WeckenTheorem}
\noindent(See \cite{ko2}, theorem 1.10).

\begin{CQ}\label{CentralQuestion}
  What happens when
  $ m \geq 2n - 2 $?
  Can we pin down extra looseness obstructions (besides the 'primary' Nielsen number
  $ (f_1, f_2) $)?
\end{CQ}

Many specific examples are known where the last claim in the Wecken theorem fails to hold as soon as the dimension assumption is not satisfied (compare e.g. the discussion of the 'Wecken condition' in \cite{ko3}, 1.18--1.29, or in \cite{kr}, table 1.18).

Already in the first critical dimension setting (when
$ m = 2n -2 $ \!\!\!)
we encounter very interesting phenomena.

\begin{MyExa}\label{example1}
  \textit{(Fixed point theory: }
  $\! f\! :=\! f_1\! $
  arbitrary selfmap of
  $ X = Y ,$
  $f_2 = \textnormal{identity map} $).
  Here
  $ m = n $,
  and in the dimension range
  $ m \geq 3 $
  our Wecken theorem implies that
  $ N(f, \id) = \MCC(f, \id) = \MC(f, \id) $
  agrees with the minimum number
  $ \MF(f) = min \left\lbrace \# C(f', \id) \MyVert f' \sim f \right\rbrace $
  of fixed points.
  This is the classical Wecken theorem (from 1941/42) for closed smooth manifolds (cf. \cite{b}, p. 12, and \cite{ko2}, pp. 225--227).
\end{MyExa}

In dimension
$ m = 2 $ \quad
J. Nielsen had already shown in the 1920s that
$ \MF(f) = N(f, \id) $
holds whenever
$ X $
is a closed connected surface with Euler characteristic
$ \chi(X) \geq 0 $.
For a long time this restriction was believed to be merely technical.

However, in 1985 B. Jiang proved that each surface
$ X $
having a strictly negative Euler characteristic allows a selfmap
$ f $
such that
$ \MF(f) \neq N(f, \id) $
(cf. \cite{j}).
Actually, later X. Zhang, M. Kelly and B. Jiang showed much more: whenever
$ \chi(X) < 0 $
the difference
$ \MF(f) - N(f, \id) $
becomes arbitrarily large for suitable selfmaps
$ f $
of the surface
$ X $
(cf. \cite{b}, p. 16).

B. Jiang used an approach via braid groups.
Could the phenomena described above also be captured by secondary obstructions? \qed
\vspace{0.7\baselineskip}
\begin{exa}\label{example2}
  $ X = S^{2n - 2},\ Y = \RP(n) $.
\end{exa}

\begin{MythmNo}
  (cf. \cite{ko3}, 1.27, or \cite{kr}, 1.13).
  Assume
  $ n $
  is even,
  $ n \neq 2, 4, 8 $.
  Let
  $ \widetilde{f} : S^{2n-2} \longrightarrow S^n $
  be a lifting of a map
  $ f : S^{2n-2} \longrightarrow \RP(n) $.
  
  Then the pair
  $ (f, f) $
  is loose if and only if both
  $ N(f,f) $
  and the Kervaire invariant
  $ K( [\widetilde{f}] ) $
  vanish.
\end{MythmNo}

Originally M. Kervaire introduced his
$\!\! \text{(}\Z_2\text{--valued)} \!\!$
invariant in order to exhibit a triangulable closed manifold which does not admit any differentiable structure (cf. \cite{ke}).
Subsequently M. Kervaire and J. Milnor used it in their classification of exotic spheres (cf. \cite{km}).
And now the Kervaire invariant makes a somewhat surprising appearance as a full secondary looseness obstruction in a very specific selfcoincidence setting. \qed
\vspace{0.7\baselineskip}

In this paper we start exploring the following general approach to constructing secondary looseness obstructions when
$ m \geq 2n-2 $.
Given a pair of maps
$ f_1, f_2 : X^m \longrightarrow Y^n $
such that
$ N(f_1, f_2) = 0 $
let us try to mimic the proof of the Wecken theorem in \cite{ko2} and measure somehow the obstacles which we encounter in the process.

We concentrate on the first critical dimension setting
$ m = 2n-2$.
Here the first two steps of the proof in \cite{ko2} (embedding a nullbordism and describing its normal bundle via a suitable desuspension) present no difficulties.
But the third step leads to a pair of maps
$ w_1, w_2 : W^n \longrightarrow Y^n $
which is coincidence free on the boundary of a given compact
$ n\text{--manifold } W $
and must be made coincidence free on all of
$ W $.
This 'secondary coincidence problem' leads to our definition of the \textit{secondary Nielsen number 
$ SecN(f_1, f_2) \!\! $},
a nonnegative integer which depends only on the homotopy classes of
$ f_1 \text{ and } f_2 $.

Clearly, if the pair
$ (f_1, f_2) $
is loose to begin with, then
\begin{equation*}
  N(f_1, f_2) = 0 = SecN(f_1, f_2) \,.
\end{equation*}
In turn we have

\begin{thm}\label{thm:theorem1}
  Assume that
  $ Y^n $
  is simply connected,
  $ n $
  even,
  $ n \neq 2 $.\\
  Given arbitrary maps
  $ f_1, f_2 : X^{2n-2} \longrightarrow Y^n $
  we have:
  the pair
  $ (f_1, f_2) $
  is loose if and only if both the ('primary') Nielsen number
  $ N(f_1, f_2) $
  and the secondary Nielsen number
  $ SecN(f_1, f_2) $
  vanish.
\end{thm}

\begin{exa}\label{example3}
  Let
  $ f : S^{2n-2} \longrightarrow S^n $
  be a map between spheres of the indicated dimensions.
  Then
  $ N(f, f) $
  is known to be a full looseness obstruction (see e.g. \cite{ko3}, 1.19).
  Hence
  $ SecN(f, f) = 0 $.
\end{exa}

\begin{exa}\label{example4}
  Let
  $ Y = \KP(n'),\ \K = \C \text{ or } \H,\, n' \geq 2 $,
  be complex or quaternionic projective space of (real) dimension
  $ n = dn' $
  where
  $ d = \dim_{\R}(\K) \in \lbrace 2, 4 \rbrace $.
  Let
  $ \partial_Y : \pi_m(Y) \longrightarrow \pi_{m-1}(S^{n-1}) $
  denote the boundary homomorphism in the exact homotopy sequence of the tangent sphere bundle
  $ \ST(Y) $
  over
  $ Y $.
  Then, given a map
  $ f : S^{2n-2} \longrightarrow Y $,
  it is well known (cf. \cite{ko3}, 1.19) that the pair
  $ (f, f) $
  is loose precisely if
  $ \partial_Y([f]) = 0 $;
  in contrast, the Nielsen number
  $ N(f,f) $
  vanishes precisely if the \textit{suspended} value
  $ E(\partial_Y([f]) \in \pi_{2n-2}(S^n) $
  is trivial.
  
  Now assume that
  $ n \neq 4, 8 $
  and that
  $ N(f, f) = 0 $.
  Then
  \begin{equation*}
    \partial_Y([f]) \in \Ker E \cong \Z_2 = \lbrace 0, 1 \rbrace
  \end{equation*}
  agrees with the secondary Nielsen number
  $ SecN(f, f) $.
  It can take a nontrivial value here if and only if
  \begin{equation*}
    \lbrace 0 \rbrace \neq \partial_Y(\pi_{2n-2}(\KP(n'))) \cap \Ker(E : \pi_{2n-3}(S^{n-1}) \longrightarrow \pi_{2n-2}(S^n))
  \end{equation*}
  (i.e. the Wecken condition for
  $ (2n-2, \KP(n')) \!\!$,
  cf. [Ko3], Definition 1.18, fails to hold), or, equivalently, if and only if
  \begin{equation*}
    0 = j_{\K^{*}}([\iota_{n-1}, \iota_{n-1}]) \in \pi_{2n-3}(V_{n'+1,2}(\K)) \,;
  \end{equation*}
  here
  \begin{equation*}
    j_{\K} \colon S^{n-1} \subset V_{n'+1,2}(\K), \ \ \ \ j_{\K}(v) := ((0, \ldots, 0, 1), (v,0)), \ \ v \in S^{n-1} \,,
  \end{equation*}
  denotes the fiber inclusion into the Stiefel manifold of orthonormal
  $ 2\!\!$--frames
  in
  $ \K^{n'+1} \!\!\!$.
\end{exa}

\section{'Primary' Nielsen numbers}

In this section we recall basic geometric facts about coincidences and the resulting definition of the Nielsen number
$ N(f_1, f_2) $
of a given pair of maps
$ f_1, f_2 \colon X^m \longrightarrow Y^n ,$
$ m, n \geq 1 \!\!$.
(For details see \cite{ko2}).

After suitable approximations we may assume that
$ (f_1, f_2) $
is a generic pair, i.e.
$ (f_1, f_2) \colon  X \longrightarrow Y \times Y $
is smooth and transverse to the diagonal
$ \Delta := \lb (y_1, y_2) \in Y \times Y \MyVert y_1 = y_2 \rb $.
Then the coincidence locus
\begin{equation*}
  C := C(f_1, f_2) = (f_1, f_2)^{-1}(\Delta)
\end{equation*}
is a closed smooth submanifold of
$ X $.
It is naturally equipped with two important geometric 'coincidence data'.
On the one hand we have a map
\begin{equation}\label{eq:equation2.1}
  \widetilde{g} \colon C \longrightarrow E(f_1, f_2) := \lb ( x, \theta) \in X \times Y^I \MyVert \theta(0) = f_1(x), \theta(1) = f_2(x) \rb
\end{equation}
such that
\begin{equation*}
  g := \pr \comp \widetilde{g} = \text{inclusion } \colon C \subset X \, ,
\end{equation*}
where
$ \pr $
denotes the projection from the 'pathspace'
$ E(f_1, f_2) $
to
$ X \!$;
$ \widetilde{g} $
is defined by
\begin{equation*}
  \widetilde{g}(x) = \lt x, \text{constant path at } f_1(x) = f_2(x) \rt, \quad x \in C\, .
\end{equation*}
On the other hand, the normal bundle
$ \nu(C, X) $
of
$ C $
in
$ X $
is described by the vector bundle isomorphism
\begin{equation}\label{eq:equation2.2Strich}
  \tag{\ref{eq:equation2.2}'}
  \nu(C, X) \xrightarrow{T(f_1, f_2)} ((f_1, f_2)\vert C)^{*} (\nu(\Delta, Y \times Y)) \cong f^{*}_1(TY) \vert C
\end{equation}
over
$ C $;
in turn this yields the (stable) isomorphism
\begin{equation}\label{eq:equation2.2}
  \overline{g} \;\colon\; TC \;\oplus\; f^{*}_1(TY) \vert C \quad \cong \quad TX \vert C \, .
\end{equation}
The resulting normal bordism class
\begin{equation*}
  \widetilde{\omega}(f_1, f_2) := \left[(C, \widetilde{g}, \overline{g}) \right] \; \in \; \Omega_{m-n}(E(f_1, f_2); \widetilde{\varphi} )
\end{equation*}
is our basic 'primary' looseness obstruction (cf. \cite{ko2}).

The decomposition of the pathspace
$ E(f_1, f_2) $
into its pathcomponents
$ A $
yields the decomposition of the coincidence manifold
$ C $
into the closed manifolds
\begin{equation*}
  C_A := \widetilde{g}^{-1}(A), \quad A \in \pi_0(E(f_1, f_2)) \, ,
\end{equation*}
(named \textit{Nielsen classes of the pair}
$ (f_1, f_2) $).
The pathcomponent
$ A \in \pi_0(E(f_1, f_2)) $
is called \textit{inessential} or \textit{essential} according as the bordism class
\begin{equation*}
  \widetilde{\omega}_A(f_1, f_2) := \left[ C_A, \widetilde{g} \vert C_A, \overline{g} \vert \right] \; \in \; \Omega_{m-n}(E(f_1, f_2); \widetilde{\varphi})
\end{equation*}
(of the coincidence data restricted to
$ C_A \!\!$)
vanishes or not.
By definition the \textit{Nielsen number}
$ N(f_1, f_2) $
is the number of essential pathcomponents
$ A $
of
$ E(f_1, f_2) \!\!\!$.
It is finite and depends only on the homotopy classes of
$ f_1 $
and
$ f_2 \!\!$.

\begin{remNo}
  The Nielsen number
  $ N(f_1, f_2) $
  is denoted by
  $ \widetilde{N}(f_1, f_2) $
  in \cite{ko3} and \cite{ko4}.
  An important role is also played by the refined (nonstabilized) version
  $ N^{\#}(f_1, f_2) $
  of our Nielsen number (cf. e.g. \cite{ko3}, \cite{ko4}).
  Furthermore, in \cite{ko4} a whole intermediate hierarchy
  \begin{equation*}
    (\MC \geq \MCC \geq ) N^{\#} \equiv N_0 \!\geq\! N_1 \!\geq\! N_2 \!\geq\! \cdots \!\geq\! N_r \!\geq\! \cdots \!\geq\! N_{\infty} \equiv \widetilde{N} \geq 0
  \end{equation*}
  of (primary) Nielsen numbers is discussed.
  However, they all coincide in the special dimension setting
  $ m = 2n - 2 $
  which will interest us in the remainder of this paper.
\end{remNo}

\section{Secondary Nielsen numbers}

Now we concentrate on the case
$ m = 2n-2\,, \; n\geq 2 $.

Thus let
$ f_1, f_2 \colon X^{2n-2} \longrightarrow Y^n $
be a generic pair.
Assume that the Nielsen coincidence class
$ C_A = \widetilde{g}^{-1}(A) $
corresponding to some pathcomponent
$ A $
of
$ E(f_1, f_2) $
is inessential.
Then we can choose a connected,
$ (n-1) \!\!$--dimensional
nullbordism
$ B $
of
$ C_A \!\!\!$,
together with maps
\begin{equation}\label{eq:equation3.1}
  \widetilde{G} \colon B \longrightarrow E(f_1, f_2) \quad \text{and} \quad
  G := \pr \comp \widetilde{G} \colon B \longrightarrow X
\end{equation}
extending
$ \widetilde{g} $
and
$ g $,
resp., (compare \eqref{eq:equation2.1}) on the one hand, as well as a (stable) vector bundle isomorphism
\begin{equation}\label{eq:equation3.2}
  \overline{G} \colon TB \oplus G^{*}(f^{*}_1(TY)) \cong G^{*}(TX) \oplus \Rschlange
\end{equation}
extending
$ \overline{g} $
(cf. \eqref{eq:equation2.2}) on the other hand.

Due to our dimension assumption the
$ n \!\!$--plane
bundle
$ G^{*}(f^{*}_1(TY)) $
allows a nowhere zero section which spans a trivial line bundle
$ \Rschlange $
over
$ B $.
In view of \eqref{eq:equation3.2} and according to Smale--Hirsch theory we can deform
$ G $
until it is an immersion.
Generically its selfintersection set consists of finitely many isolated points.
We may 'push' them all along suitable arcs across the boundary
$ \partial B = C $.
So in the end
$ G $
is a smooth embedding which extends the inclusion
$ C \subset X $.

Next let
$ \nu(G) $
be the normal bundle of
$ G $.
Compose the obvious isomorphism
\begin{equation*}
  TB \oplus \nu(G) \oplus \Rschlange \ \cong \ TX \vert B \oplus \Rschlange
\end{equation*}
with the isomorphism
$ \overline{G} $
in \eqref{eq:equation3.2}.
Again in view of our dimension assumption we can desuspend to get the isomorphism
\begin{equation}\label{eq:equation3.2Strich}
  \tag{\ref{eq:equation3.2}'}
  \nu(G) \oplus \Rschlange \cong f^{*}_1(TY) \vert B
\end{equation}
which extends \eqref{eq:equation2.2Strich} when we put
$ \nu(C_A, B) = \Rschlange \vert C_A \!$.

As in the proof of the Wecken therorem in \cite{ko2}, pp. 223--224, let
$ i_2 \colon B \longrightarrow [0, 1] $
be a smooth function which is essentially defined by the second projection on a collar
$ C_A \times [0, \frac{1}{2}] $
of
$ C_A $
in
$ B $,
and takes the constant value
$ \frac{1}{2} $
outside of this collar.
Consider the embedding
\begin{equation*}
  i := (G, i_2) \colon B \hookrightarrow X \times [0, 1]
\end{equation*}
which extends the embedding
$ g \colon \partial B = C_A \hookrightarrow X = X \times \lb 0 \rb $
(compare \cite{ko1}, figure 3.8).
As in \cite{ko2}, pp. 223--224, the isomorphism in \eqref{eq:equation3.2Strich} allows us to extend
$ (f_1, f_2) $,
defined on
$ X = X \times \lb 0 \rb \!\!$,
to a map
\begin{equation*}
  (F_1, F_2) \colon X \times \lb 0 \rb \cup T \longrightarrow Y \times Y
\end{equation*}
where
$ T $
is a suitable (compact) tubular neighborhood of
$ i(B) $
in
$ X \times [0,1] $;
the coincidence locus of the pair
$ (F_1, F_2) $
is just
$ i(B) $.

Now let
$ W \subset X \times I $
denote that part of the 'shadow'
\begin{equation*}
  R_1 = \lb (x, t) \in G(B) \times I \MyVert 0 \leq t \leq i_2(x) \rb
\end{equation*}
(compare \cite{ko1}, pp. 38--39) which lies outside of the interior
$ \mathring{T} $
of
$ T $.
The map
$ \widetilde{G} $
(cf. \eqref{eq:equation3.1}) determines an extension of the coincidence free pair
$ (F_1 \vert \partial W, F_2 \vert \partial W) $
to the whole
$ n \!\!$--dimensional
manifold
$ W $.
But now coincidences (which are generically
$ 0 \!\!$--dimensional)
may occur in the interior of
$ W $.
They lead to
$ A \in \pi_0(E(f_1, f_2)) $
being called either \textit{
$ 2 \!\!$--essential}
(i.e. \textit{'essential of second order'}) or not.
If 
$ C_A $
can be made empty by suitable homotopies of
$ f_1 $
and
$ f_2 $,
then
$ A $
is certainly
$ 2 \!\!$--inessential.

Now assume that the (primary) Nielsen number
$ N(f_1, f_2) $
vanishes.
Then each pathcomponent
$ A $
of
$ E(f_1, f_2) $
is inessential in the classical (primary) sense, and hence either
$ 2 \!\!$--essential
or not.
\textit{We define the secondary Nielsen number
$ SecN(f_1, f_2) $
to be the number of
$ 2 \!\!$--essential
pathcomponents
$ A \in \pi_0(E(f_1, f_2)) $.
}
Clearly 
$ SecN(f_1, f_2) $
is a nonnegative integer smaller or equal to the (not necessarily finite) number of pathcomponents of
$ E(f_1, f_2) $
(which is also known as the
\textit{Reidemeister number of
$ (f_1, f_2) \!\!\!$;
}
see \cite{ko2}, 2.1, for explicit calculations).

If
$ \pi_1(Y) = 0 $,
then not only
$ E(f_1, f_2) $
but also
$ E(F_1, F_2) $
(in the construction discussed above) is pathconnected.
This allows us to prove the theorem \ref{thm:theorem1} of the introduction.
Details and generalizations will be given in a future paper.


\end{document}